\documentclass[11pt]{article}
\usepackage{amsfonts}
\usepackage{amsfonts,mathrsfs,amssymb,amsthm,mathptm}
\usepackage{amsmath,amscd}
\usepackage{mathptm,pslatex}
\usepackage{color}
\usepackage[dvips]{graphicx}
\usepackage[all]{xy}
\usepackage{graphicx}
\usepackage{fancyhdr}

\begin{document}
\newtheorem{Def}{Definition}[section]
\newtheorem{Bsp}[Def]{Example}
\newtheorem{Prop}[Def]{Proposition}
\newtheorem{Theo}[Def]{Theorem}
\newtheorem{Lem}[Def]{Lemma}
\newtheorem{Koro}[Def]{Corollary}
\newtheorem{Nota}[Def]{Notation}
\theoremstyle{definition}
\newtheorem{Rem}[Def]{Remark}

\newcommand{\add}{{\rm add}}
\newcommand{\con}{{\rm con}}
\newcommand{\gd}{{\rm gldim}}
\newcommand{\repdim}{{\rm repdim}}
\newcommand{\sd}{{\rm stdim}}
\newcommand{\sr}{{\rm sr}}
\newcommand{\dm}{{\rm domdim}}
\newcommand{\cdm}{{\rm codomdim}}
\newcommand{\tdim}{{\rm dim}}
\newcommand{\E}{{\rm E}}
\newcommand{\K}{{\rm k}}
\newcommand{\Mor}{{\rm Morph}}
\newcommand{\End}{{\rm End}}
\newcommand{\ind}{{\rm ind}}
\newcommand{\rsd}{{\rm res.dim}}
\newcommand{\rd} {{\rm rd}}
\newcommand{\ol}{\overline}
\newcommand{\overpr}{$\hfill\square$}
\newcommand{\rad}{{\rm rad}}
\newcommand{\soc}{{\rm soc}}
\renewcommand{\top}{{\rm top}}
\newcommand{\pd}{{\rm pdim}}
\newcommand{\id}{{\rm idim}}
\newcommand{\fld}{{\rm fdim}}
\newcommand{\Fac}{{\rm Fac}}
\newcommand{\Gen}{{\rm Gen}}
\newcommand{\fd} {{\rm fin.dim}}
\newcommand{\Fd} {{\rm Fin.dim}}
\newcommand{\Pf}[1]{{\mathcal P}^{<\infty}(#1)}
\newcommand{\DTr}{{\rm DTr}}
\newcommand{\cpx}[1]{#1^{\bullet}}
\newcommand{\D}[1]{{\mathcal D}(#1)}
\newcommand{\Dz}[1]{{\mathcal D}^+(#1)}
\newcommand{\Df}[1]{{\mathcal D}^-(#1)}
\newcommand{\Db}[1]{{\mathcal D}^b(#1)}
\newcommand{\C}[1]{{\mathcal C}(#1)}
\newcommand{\Cz}[1]{{\mathcal C}^+(#1)}
\newcommand{\Cf}[1]{{\mathcal C}^-(#1)}
\newcommand{\Cb}[1]{{\mathcal C}^b(#1)}
\newcommand{\Dc}[1]{{\mathcal D}^c(#1)}
\newcommand{\Kz}[1]{{\mathcal K}^+(#1)}
\newcommand{\Kf}[1]{{\mathcal  K}^-(#1)}
\newcommand{\Kb}[1]{{\mathcal K}^b(#1)}
\newcommand{\DF}[1]{{\mathcal D}_F(#1)}
\newcommand{\rpd}{{\rm repdim}}

\newcommand{\Kac}[1]{{\mathcal K}_{\rm ac}(#1)}
\newcommand{\Keac}[1]{{\mathcal K}_{\mbox{\rm e-ac}}(#1)}

\newcommand{\modcat}{\ensuremath{\mbox{{\rm -mod}}}}
\newcommand{\Modcat}{\ensuremath{\mbox{{\rm -Mod}}}}
\newcommand{\Spec}{{\rm Spec}}

\newcommand{\stmc}[1]{#1\mbox{{\rm -{\underline{mod}}}}}
\newcommand{\Stmc}[1]{#1\mbox{{\rm -{\underline{Mod}}}}}
\newcommand{\prj}[1]{#1\mbox{{\rm -proj}}}
\newcommand{\inj}[1]{#1\mbox{{\rm -inj}}}
\newcommand{\Prj}[1]{#1\mbox{{\rm -Proj}}}
\newcommand{\Inj}[1]{#1\mbox{{\rm -Inj}}}
\newcommand{\PI}[1]{#1\mbox{{\rm -Prinj}}}
\newcommand{\GP}[1]{#1\mbox{{\rm -GProj}}}
\newcommand{\GI}[1]{#1\mbox{{\rm -GInj}}}
\newcommand{\gp}[1]{#1\mbox{{\rm -Gproj}}}
\newcommand{\gi}[1]{#1\mbox{{\rm -Ginj}}}

\newcommand{\opp}{^{\rm op}}
\newcommand{\otimesL}{\otimes^{\rm\mathbb L}}
\newcommand{\rHom}{{\rm\mathbb R}{\rm Hom}\,}
\newcommand{\pdim}{\pd}
\newcommand{\Hom}{{\rm Hom}}
\newcommand{\Coker}{{\rm Coker}}
\newcommand{ \Ker  }{{\rm Ker}}
\newcommand{ \Cone }{{\rm Con}}
\newcommand{ \Img  }{{\rm Im}}
\newcommand{\Ext}{{\rm Ext}}
\newcommand{\StHom}{{\rm \underline{Hom}}}
\newcommand{\StEnd}{{\rm \underline{End}}}

\newcommand{\KK}{I\!\!K}

\newcommand{\gm}{{\rm _{\Gamma_M}}}
\newcommand{\gmr}{{\rm _{\Gamma_M^R}}}

\def\vez{\varepsilon}\def\bz{\bigoplus}  \def\sz {\oplus}
\def\epa{\xrightarrow} \def\inja{\hookrightarrow}

\newcommand{\lra}{\longrightarrow}
\newcommand{\llra}{\longleftarrow}
\newcommand{\lraf}[1]{\stackrel{#1}{\lra}}
\newcommand{\llaf}[1]{\stackrel{#1}{\llra}}
\newcommand{\ra}{\rightarrow}
\newcommand{\dk}{{\rm dim_{_{k}}}}

\newcommand{\holim}{{\rm Holim}}
\newcommand{\hocolim}{{\rm Hocolim}}
\newcommand{\colim}{{\rm colim\, }}
\newcommand{\limt}{{\rm lim\, }}
\newcommand{\Add}{{\rm Add }}
\newcommand{\Prod}{{\rm Prod }}
\newcommand{\Tor}{{\rm Tor}}
\newcommand{\Cogen}{{\rm Cogen}}
\newcommand{\Tria}{{\rm Tria}}
\newcommand{\Loc}{{\rm Loc}}
\newcommand{\Coloc}{{\rm Coloc}}
\newcommand{\tria}{{\rm tria}}
\newcommand{\Con}{{\rm Con}}
\newcommand{\Thick}{{\rm Thick}}
\newcommand{\thick}{{\rm thick}}
\newcommand{\Sum}{{\rm Sum}}

{\Large \bf
\begin{center}
Structure of Terwilliger algebras of quasi-thin association schemes
\end{center}}

\medskip
\centerline{\textbf{Zhenxian Chen}  and \textbf{Changchang Xi}$^*$ }

 \renewcommand{\thefootnote}{\alph{footnote}}
 \setcounter{footnote}{-1} \footnote{ $^*$ Corresponding author.
 Email: xicc@cnu.edu.cn; Fax: 0086 10 68903637.}
 \renewcommand{\thefootnote}{\alph{footnote}}
 \setcounter{footnote}{-1} 
 \renewcommand{\thefootnote}{\alph{footnote}}
 \setcounter{footnote}{-1} 

 \begin{abstract}
We show that the Terwilliger algebra of a quasi-thin association scheme over a field is always a quasi-hereditary cellular algebra in the sense of Cline-Parshall-Scott and of Graham-Lehrer, repsectively, and that the basic algebra of the Terwilliger algebra is the dual extension of a star with all arrows pointing to its center if the field has characteristic $2$. Thus many homological and representation-theoretic properties of these Terwilliger algebras can be determined completely.
For example, the Nakayama conjecture holds for Terwilliger algebras of quasi-thin association schemes.
 \end{abstract}


\section{Introduction}
Association schemes have played an important role in the theory of algebraic combinatorics and designs (see \cite{et, phz}). To understand them algebraically, Terwilliger associated each association scheme with an algebra over a field or more generally, over a commutative ring (see \cite{tw1,tw2,tw3}). This algebra nowadays is termed as the Terwilliger algebra (for example, see \cite{phz}). There is a large variety of literature on the study of Terwilliger algebras, especially, for special schemes or ground fields (for example, see \cite{hm}, \cite{lm}, \cite{mp}, \cite{lmw} and the references therein). Recently, Hanaki studies the modular case of Terwilliger algebras \cite{ah}, and Jiang investigates the Jacobson radicals of the Terwilliger algebras of quasi-thin association schemes \cite{yj1}. The tools used in their study are mainly combination of combinatorics with ring theory.

The purpose of this note is to understand the Terwilliger algebras of quasi-thin association schemes from the view point of representation theory of algebras, namely we show the following structural and homological result.

\begin{Theo} \label{thm}Let $R$ be a field of arbitrary characteristic and $S$ a quasi-thin association scheme on a finite set. Then the Terwilliger $R$-algebra of $S$ is quasi-hereditary in the sense of Cline-Parshall-Scott, and cellular in the sense of Graham-Lehrer. Moreover, if the field $R$ has characteristic $2$, the basic algebra of the Terwilliger $R$-algebra of $S$ is the dual extension of a star, and has global dimension at most $2$.
\end{Theo}

For a finite-dimensional algebra $A$ over a field, let $P_1, \cdots, P_n$ be a complete set of non-isomorphic, indecomposable projective $A$-modules. Then the basic algebra $\Lambda$ of $A$ is defined to be the endomorphism algebra of the $A$-module $\bigoplus _{i=1}^n P_i$. It is known that $A$ and $\Lambda$ are Morita equivalent, that is, they have the equivalent module categories.

Observe that the notion of cellular algebras in the sense of Graham-Lehrer in \cite{gl} is completely different from the one in \cite[Section D, p.23]{weis}.

The proof of Theorem \ref{thm} is given in the next section.

\section{Proof of the result}\label{Preliminaries}
In this section we first recall basics on Terwilliger algebras of schemes and then prove the main result. Also, representation-theoretic and homological properties of Terwilliger algebras of quasi-thin association schemes are deduced.

\subsection{Terwilliger algebras of association schemes}\label{SI2.2}

Throughout this note, $R$ is a field. For $n\in \mathbb{N}$, the symbol $[n]$ denotes the  set $\{0,1,2, \cdots,n\}\subseteq \mathbb{N}$. The cardinality of a set $X$ is denoted by $|X|$.

\begin{Def}\label{Scheme}
An \emph{association scheme} or simply a \emph{scheme} of size $d$ on a nonempty finite set $X$ is a partition $S=\{R_0,R_1,\ldots,R_d\}$ of the Cartesian product  $X\times X$ with all parts $R_i$ nonempty, satisfying the conditions

        $(S1)$ $R_0=\{(x,x)\mid x\in X\}$.

        $(S2)$ For $i\in[d]$, there exists $i^\prime\in [d]$ such that $R_{i^{\prime}}=\{(x,y)\mid (y,x)\in R_i\}$, and

        $(S3)$ For $i,j,\ell\in [d]$ and $(x,y),(u,v)\in R_\ell$, the following holds:

        $$\left|\{z \in X \mid (x,z)\in R_i, (z,y)\in R_j\}\right|=\left|\{w\in X\mid (u,w)\in R_i,(w,v)\in R_j\}\right|.$$
\end{Def}

Now, let $S=\{R_0,R_1,\ldots,R_d\}$ be a scheme of size $d$. An element $x\in X$ is called a \emph{vertex} of $S$, and the part $R_i$ is called a \emph{relation} of $S$. For $i,j,\ell\in [d]$ and $(x,y)\in R_\ell$, we define
$$p_{ij}^\ell:=\left|\{z\in X\mid (x,z)\in R_i,(z,y)\in R_j\}\right| \in \mathbb{N}.$$

It is known from $(S3)$ that $p_{i j}^{\ell}$ is independent on the choice of $(x,y)$ in $R_{\ell}$. The number $p_{i j}^{\ell}$ is called the \emph{intersection number} of $S$ with respect to the triple $(R_i, R_j, R_\ell)$. By Definition \ref{Scheme}(S2), the $i'$ is uniquely determined by $i$. So, the \emph{valency} $k_i$ of $R_i$ is defined by $k_i:= p_{i i^{\prime}}^0$. Note that $k_i$ is just the cardinality of the set $x R_i:=\left\{y \in X\mid (x, y) \in R_i\right\}$ for any $x\in X$. Thus $k_i>0$ for $i\in [d]$ and $k_0=1$. For $j\ge 1$, let $\mathcal{A}_j:=\{i\in[d]\mid k_i=j\}$. Then $\mathcal{A}_1\neq \emptyset$.

\begin{Def} A scheme $S$ of size $d$ is called a \emph{thin scheme} if $k_i=1$ for all $i\in [d]$; and a \emph{quasi-thin} scheme if $k_i\le 2$ for all $i\in [d]$.
\end{Def}

Quasi-thin schemes was introduced in \cite{hm}, but the first result on quasi-thin schemes goes back to \cite{weis}, where it was proved that any primitive quasi-thin scheme is Schurian.

By definition, if $S$ is thin, then $\mathcal{A}_j=\emptyset$  for $ 2\le j\in \mathbb{N}$.
The following properties of schemes are well known.

\begin{Lem}{\rm \cite{dgh}}
        Let $S$ be a scheme of size $d$. If $i,j,\ell\in [d]$, then
        \begin{enumerate}
        $(1)$ $p^\ell_{ji}=p^{\ell^{\prime}}_{i^{\prime}j^{\prime}}$ and $k_i=k_{i^{\prime}}$.

        $(2)$ $\sum_{\ell=0}^{d} p^j_{i\ell}=k_i$.

        $(3)$ $k_ik_j=\sum_{\ell=0}^{d} p^\ell_{ij}k_\ell$.

        $(4)$ $k_\ell p^\ell_{ij}=k_ip^i_{\ell j^{\prime}}=k_jp^j_{i^{\prime}\ell}$.
        \end{enumerate}
    \end{Lem}

For any nonempty subsets $U, V$ of $S$, the multiplication of $U$ and $V$ is defined by $$UV:=\left\{R_\ell \in S \mid \exists\, R_i \in U, \exists\; R_j \in V, \mbox{ such that } p_{ij}^\ell>0\right\}.$$
For $i,j\in[d]$, $\left|R_i R_j\right|\le \gcd(k_i,k_j)$ for $i,j\in [d]$ by \cite[Lemma 1.5.2]{phz}, where $\gcd(m,n)$ means the greatest common divisor of $m$ and $n$.

\medskip
Let $M_X(R)$ be the $X\times X$ matrix algebra over $R$. We denote by $I$ the identity matrix in $M_X(R)$, $E_{xy}$ the matrix units for $x,y\in X$, and $J$ the matrix with all entries equal to $1$.

For a part $R_i$ of a scheme $S$, there is associated an \emph{adjacency matrix} $A_i:=(a_{xy})\in M_X(R)$ defined by $a_{xy}=1$ if $(x, y) \in R_i$, and $0$ otherwise. Thus $A_0 = I, A_i^t = A_{i^{\prime}}$ and $\sum_{i=0}^d A_i=J$, where $A^t$ is the transpose matrix of $A$. Moreover, for $i, j \in [d]$, $$A_i A_j=\sum_{\ell=0}^d p_{i j}^{\ell} A_{\ell}.$$

Fix an $x \in X$ and $i\in [d]$, the \emph{dual idempotent} of $R_i$ with respect to $x$ is defined by $E_i^*(x):=\sum_{y \in x R_i} E_{yy}\in M_X(R)$. Then

$$(\dag)\quad E_i^*(x) E_j^*(x)=\delta_{i j} E_j^*(x), \; \sum_{i=0}^d E_i^*(x)=I, \; \mbox{ and } J E_i^*(x) J = k_i J,$$where $\delta_{ij}$ is the Kronecker symbol. Thus $E^*_0, E^*_1, E^*_2, \cdots, E^*_d$ form a complete set of pairwise orthogonal idempotent elements of $M_X(R)$. Moreover, for $M=(m_{xy}) \in M_X(R)$,
$$(\ddag)\quad E_i^*(x) M E_j^*(x)=\sum_{y \in x R_i} \; \sum_{z \in x R_j} m_{y z} E_{y z}.$$
This shows $E_i^*(x) A_\ell E_j^*(x)$ is a $(0,1)$-matrix for $i,j,\ell\in[d]$.

\begin{Def}{\rm \cite{yj1}}\label{bp}
Let $n\in\mathbb{N}$. If there are numbers $i_b,j_b,\ell_b\in [d]$ for all $b\in[n]$, such that

$(1)$ $k_{i_b}=k_{\ell_b}=2$ and $p_{i_bj_b}^{\ell_b}=1$ for all $b\in [n]$, and

$(2)$ $\left| R_{i_0^\prime}R_{\ell_n} \right|=1$ and $\ell_c=i_{c+1}$ for all $c\in[n-1]$,

\noindent then the pair $(i_0,\ell_n)$ is called a bad pair of $S$.
\end{Def}

Let $\mathcal{S}$ denote the set of all bad pairs of $S$, $\mathcal{R}:=\left\{(i, j)\in \mathcal{A}_2\times \mathcal{A}_2\mid \left|R_{i^{\prime}} R_j\right|=2\right\}$, and $\mathcal{U}:=\mathcal{R} \cup \mathcal{S}$.

One defines a relation on $\mathcal{A}_2$: For $i,j\in \mathcal{A}_2$, $i\sim j$ if and only if $(i,j)\in \mathcal{U}$. Then $\sim$ is an equivalence relation on $\mathcal{A}_2$ by \cite[Lemma 7.2]{yj1}. The set of equivalence classes of $\sim$ is denoted by $\{\mathcal{C}_1, \ldots, \mathcal{C}_r\}$ for $r$ a natural number. We define $\mathcal{C}_0:=[d]$ and $b_{i j}^0(x):=E_i^*(x) J E_j^*(x)$ for $i,j\in \mathcal{C}_0$ and $x\in X$.

Now, we choose a fixed $x\in X$, a total order $\prec$ for $X$, and take $i, j \in \mathcal{C}_\ell$ for $1\le \ell\le r$. If $x R_i=\left\{y_1, y_2\right\}$ and $x R_j=\left\{z_1, z_2\right\}$ such that $y_1 \prec y_2$ and $z_1 \prec z_2$, then we define $b_{i j}^\ell(x):=E_{y_1 z_1}+E_{y_2 z_2}$. Clearly, $(b_{i j}^\ell(x))^t=b_{ji}^\ell(x)$ for $\ell\in [r]$, and $b_{j j}^{\ell}(x)=E^*_j(x)$ for $\ell=0$ and $j\in \mathcal{A}_1$, or $1\le\ell\le r$ and $j\in \mathcal{C}_{\ell}$. Let $\mathcal{B}_\ell(x):=\{b_{i j} ^ \ell(x) \mid i, j\in \mathcal{C}_\ell\}$ for $0\le \ell\le r$ and $\mathcal{B}(x):=\bigcup_{\ell=0}^r\mathcal{B}_\ell(x)$.

\begin{Def}{\rm \cite{tw1}} Let $S$ be a scheme on a finite set $X$, and $R$ be a commutative ring $R$ with identity. The Terwilliger $R$-algebra $\mathcal{T}_{R}(x)$ of a scheme $S$ on $X$ with respect to $x\in X$ is the $R$-subalgebra of $M_X(R)$ generated by $A_0, A_1, \ldots, A_d$; $E_0^*(x), E_1^*(x), \ldots,$ and $E_d^*(x)$.
\end{Def}

By a Terwilliger algebra of $S$, we means the Terwilliger algebra $\mathcal{T}_{R}(x)$ with respect to $x\in X$. Terwilliger algebras were first introduced and studied by P. Terwilliger in a series of papers \cite{tw1,tw2,tw3} for commutative schemes under a different name. These algebras now are called \emph{Terwilliger algebras}. They were studied in \cite{ah} over a field of positive characteristic under the name \emph{modular Terwilliger algebras}.

Clearly, the transpose of matrices is an $R$-linear involution of $\mathcal{T}_{R}(x)$, and  $\mathcal{T}_{R}(x)= M_X(R)$ if $S$ is thin. In general, for $x,y\in X$ with $x\ne y$, we do not have  $\mathcal{T}_{R}(x)\simeq \mathcal{T}_{R}(y)$ as algebras \cite[Section 5.1]{ah}. However, this holds true for quasi-thin scheme,  namely if $S$ is a quasi-thin scheme, then
$\mathcal{T}_R(x)\simeq \mathcal{T}_R(y)$ as $R$-algebras for $x,y \in X$ by \cite[Theorem D]{yj1}.

\textbf{From now on we assume that $S$ is a quasi-thin scheme} of size $d$ on a finite set $X$. We fix an $x\in X$ and write $\mathcal{T}:=\mathcal{T}_R(x)$ and $E_i^*:=E_i^*(x)$ for $i\in [d]$. Similarly, we write other notation involving $x$.

\begin{Prop} {\rm \cite{yj1}} \label{yj1}
Let $R$ be a field of characteristic $p\ge 0$ and $S$ be a quasi-thin scheme. Then

$(1)$ $\mathcal{B}$ is an $R$-basis of $\mathcal{T}$. Thus $\dim_R(\mathcal{T})= |\mathcal{R}|+|\mathcal{S}|+(d+1)^2$.

$(2)$  For $b_{ij}^\ell\in\mathcal{B}_{\ell}, b_{uv}^w\in \mathcal{B}_w$ with $\ell, w \in [d],$ the following holds

$$b_{ij}^\ell b_{uv}^w=\begin{cases} \delta_{ju}k_j b_{iv}^0, & \text { if } \ell=w=0, \\ \delta_{ju}\delta_{\ell w}b_{iv}^\ell, &\text { if } \ell\neq 0 \text{ and } w\neq 0, \\ \delta_{ju}b_{iv}^0, &\text { otherwise.}   \end{cases}$$

$(3)$ $ \mathcal{T} $ is semisimple if and only if $ p \neq 2 $ or $ p = 2 $ and $ S $ is thin.

$(4)$ If $p = 2$, the Jacobson radical of $\mathcal{T}$, denoted by $\rad(\mathcal{T})$, is spanned $R$-linearly by $\{b_{ij}^0\mid i,j\in [d],\max \left\{k_i, k_j\right\}=2\}$.
\end{Prop}

\subsection{Cellular  and quasi-hereditary algebras}

Let us recall the definition of cellular algebras in \cite{gl}.
\begin{Def}{\rm \cite{gl}}
Let $R$ be a  commutative noetherian domain. A unitary $R$-algebra $A$ is called a \emph{cellular algebra} with cell datum $(I,M,C,t)$ if the following are satisfied:

$(C1)$ The finite set $I$ is partially ordered. Associated with each $\lambda\in I$ there is a finite set $M(\lambda)$. The algebra $A$ has an $R$-basis $C_{S,T}^\lambda$ where $(S,T)$ run through all elements of $M(\lambda)\times M(\lambda)$ for all $\lambda\in I$.

$(C2)$ The map $t$ is an $R$-involution on $A$, that is, an $R$-linear anti-automorphism on $A$ of oder $2$, such that $C_{S,T}^\lambda$ is sent to $C_{T,S}^\lambda$.

$(C3)$ For $a\in A$, there holds $$aC_{S,T}^\lambda= \sum_{U\in M(\lambda)}r_a(U,S)C_{U,T}^\lambda+r^\prime$$where $r_a(U,S)\in R$ does not depend on $T$ and $r^\prime$ is a linear combination of basis elements $C^{\mu}_{W,V}$ with $\mu<\lambda$.
\label{def-ca}
\end{Def}

A trivial example of cellular algebras is $M_n(R)$ with $I:=\{*\}$, $M(*):=[n], C^*_{i j}:=E_{ij}$ and $t$ being the matrix transpose. Cellular algebras capture many interesting classes of algebras such as Brauer algebras \cite{gl}, Hecke algebras of
finite type \cite{geck}, centralizer matrix algebras \cite{xz2021}.

The following is an equivalent definition of cellular algebras:

\begin{Def}{\rm \cite{sx1}}\label{ca2}
Let $A$ be an $R$-algebra with with $R$ a commutative noetherian domain. Assume that there is an involution $t$ on $A$. A two-sided ideal $J$ in $A$ is called a \emph{cell ideal} if and only if $(J)t=J$ and there exists a left ideal $\Delta\subset J$ of $A$ such that $\Delta$ is finitely generated and free over $R$ and that there is an isomorphism of $A$-bimodules $\alpha:J\simeq \Delta\otimes_R(\Delta)t$ $($where $(\Delta)t\subset J$ is the image of $\Delta$ under $t$) making the following diagram commutative:

$$\xymatrix{
J \ar[r]^{\rlap{$\scriptstyle \alpha$}\hspace{2em}} \ar[d]_t & \Delta\otimes_R(\Delta)t \ar[d]^{x\otimes y\mapsto(y)t\otimes(x)t} \\
J \ar[r]_{\rlap{$\scriptstyle \alpha$}\hspace{2em}} & \Delta\otimes_R(\Delta)t
}$$

The algebra $A$ $($with the involution $t$$)$ is called \emph{cellular} if and only if there is an $R$-module decomposition $A = J_1^\prime\oplus J_2^\prime\oplus \cdots\oplus J_n^\prime$ $($for some $n$$)$ with $(J_j^\prime)t = J_j^\prime$ for each $j$, such that setting $J_j=\oplus_{l=1}^j J_l^\prime$ gives a chain of two-sided ideals of $A:0=J_0\subset J_1\subset J_2\subset\cdots\subset J_n = A$ $($each of them fixed by $t$$)$ and that, for each $1\le j\le n$, the quotient $J_j^\prime=J_j /J_{j-1}$ is a cell ideal $($with respect to the involution induced by $t$ on the quotient$)$ of $A/J_{j-1}$.
    \end{Def}

Next, we recall the definition of quasi-hereditary algebras introduced in \cite{cps}.

\begin{Def}{\rm \cite{cps}}\label{qha}
Let $A$ be a finite-dimensional algebra over a field. An ideal $J$ in $A$ is called a \emph{heredity ideal} if $J$ is idempotent, $J(\rad(A))J = 0$ and $J$ is a projective left (or right) $A$-module. The algebra $A$ is said to be \emph{quasi-hereditary} provided there is a finite chain $0=J_0\subset J_1\subset J_2\subset\cdots\subset J_n = A$ of ideals in $A$ such that $J_j/J_{j-1}$ is a heredity ideal in $A/J_{j-1}$ for all $j$. Such a chain is then called a \emph{heredity chain} of the quasi-hereditary algebra $A$.
\end{Def}
To judge whether a cellular algebra is quasi-hereditary, we have the following criterions.

\begin{Lem}{\rm \cite[Lemma 2.1]{sx2}}\label{ciqh}
Let $A$ be a finite-dimensional cellular algebra over a field, with a cell chain $0=J_0\subset J_1\subset\cdots \subset J_n=A$. Then the given cell chain is a heredity chain $($make $A$ into a quasi-hereditary algebra$)$ if and only if all $J_l$ satisfy $J_l^2\nsubseteq J_{l-1}$ if and only if $n$ equals the number of isomorphism classes of simple modules.
\end{Lem}

\subsection{Proof of the main result}

This section is devoted to the proof of the following main result.

\begin{Theo}\label{tiqhc} Let $R$ be a field of characteristic $p\ge 0$, and assume that $S$ be a quasi-thin scheme on a finite set $X$. Then the following hold.

$(1)$ The Terwilliger $R$-algebra $\mathcal{T}$ of $S$ is a quasi-hereditary, cellular algebra with respect to the matrix transpose.

$(2)$ If $R$ is of characteristic $2$ and $\mathcal{A}_2$ has $r\ge 0$ equivalence classes, then the basic algebra $\Gamma$ of Terwilliger $R$-algebra $\mathcal{T}$ of $S$ is isomorphic to an $R$-algebra $\Lambda$ given by the following quiver
$$
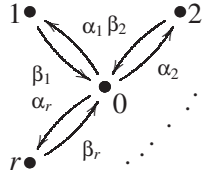
\xy{
(12,-1)*+{\cdot};(11,-2.8)*+{\cdot};(8.5,-5.2)*+{\cdot};(6.2,-7.5)*+{\cdot};(4.9,-9)*+{\cdot};
(3,-10)*+{\cdot};(2,-2)*+{0};(-12,10)*+{1};(12,10)*+{2};(-12,-10)*+{r};
\ar@/_0.3pc/_{\beta_1} (0,0)*+{\bullet}; (-10,10)*+{\bullet};
\ar@/_0.3pc/_{\alpha_1} (-10,10)*+{\bullet}; (0,0)*+{\bullet};
\ar@/_0.3pc/_{\beta_2} (0,0)*+{\bullet}; (10,10)*+{\bullet};
\ar@/_0.3pc/_{\alpha_2} (10,10)*+{\bullet}; (0,0)*+{\bullet};
\ar@/_0.3pc/_{\beta_r} (0,0)*+{\bullet}; (-10,-10)*+{\bullet};
\ar@/_0.3pc/_{\alpha_r} (-10,-10)*+{\bullet}; (0,0)*+{\bullet};
}\endxy$$
with relations $\alpha_i\beta_i=0$ for $1\le i\le r$.
\end{Theo}

{\it Proof.} (1) We apply Proposition \ref{yj1} to show the cellularity of the Terwilliger $R$-algebra $\mathcal{T}$ of a quasi-thin scheme $S$.

If $\mathcal{A}_2=\emptyset$, that is $r=0$, then $\mathcal{T}= M_X(R)$ with $X$ a finite set. Thus Theorem \ref{tiqhc} follows.

Now, assume $\mathcal{A}_2\neq\emptyset$. We define  $I:=[r]$ with the partial order $\{0\prec 1,0\prec 2,\cdots, 0\prec r\}$, $M(\ell):=\mathcal{C}_\ell$ and $t$ being the transpose of matrices, namely $a\mapsto a^t,\forall a\in \mathcal{T}$. We show that $(I,M, b,t)$ is a cell datum for $\mathcal{T}$.

Indeed, by Proposition \ref{yj1}, $\{b_{i j} ^\ell \mid \ell\in I, i, j\in M(\ell)\}$ is an $R$-basis of $\mathcal{T}$. Clearly, $t$ is an $R$-linear anti-automorphism of order $2$ on $\mathcal{T}$ and $(b^{\ell}_{ij})^t=b^{\ell}_{ji}$ for $\ell\in I$ and $i,j\in M(\ell)$. It remains to consider the multiplication of basis elements of $\mathcal{T}$ and verify Definition \ref{def-ca}(C3). But this follows immediately from Proposition \ref{yj1}.

To see that $\mathcal{T}$ is quasi-hereditary, we define a chain of ideals in $\mathcal{T}$.
Let $J_0=0$ and $J_i$ be the $R$-span of the elements in $\bigcup_{\ell=0}^{i-1}\mathcal{B}_\ell$ for $0\neq i\in [r+1]$. Thus the chain $0=J_0\subset J_1\subset \cdots \subset J_{r+1}= \mathcal{T}$ is a cell chain of ideals in $\mathcal{T}$ (see \cite{sx1}). To show that this cell chain is a heredity chain, we prove the following:

(a) For $0\in \mathcal{A}_1$, we get $b_{00}^0\in \mathcal{B}_0\subset J_1$ and $b_{00}^0b_{00}^0=k_0b_{00}^0=b_{00}^0\neq 0$.

(b) For $\ell\in \{1,2, \cdots, r\}$, it follows from $\mathcal{C}_\ell\neq \emptyset$ that there is an element $u_\ell\in \mathcal{C}_\ell$, such that $b_{u_\ell u_\ell}^\ell\in \mathcal{B}_\ell\subset J_{\ell+1}$, and $$b_{u_\ell u_\ell}^\ell b_{u_\ell u_\ell}^\ell=b_{u_\ell u_\ell}^\ell\notin J_\ell.$$
Thus the chain $0=J_0\subset J_1\subset \cdots \subset J_{r+1}= \mathcal{T}$ is a heredity chain by Lemma \ref{ciqh}, and therefore $\mathcal{T}$ is a quasi-hereditary algebra.

(2) If $\mathcal{T}$ is semisimple, then $r=0$ and Theorem \ref{tiqhc} is clearly true.  Now suppose that $\mathcal{T}$ is not semisimple, that is, we are in the case (2) of Theorem \ref{tiqhc} with $r>0$.

For $\ell\in [r]$, let  $\mathcal{D}_0:=\mathcal{A}_1, \mathcal{D}_\ell:=\mathcal{C}_\ell$ for $\ell\neq 0$. We have the following

(i) $E:=\{b_{ii}^\ell\mid \ell\in [r], i\in \mathcal{D}_\ell\}$ is a complete set of pairwise orthogonal idempotents of $\mathcal{T}$ by $(\dag)$.

(ii) For $\ell\in [r],i\in \mathcal{D}_\ell$, let $P^{\ell}_i:=\mathcal{T}b^{\ell}_{ii}$. Then $P^{\ell}_i\simeq P_j^{\ell}$ as $\mathcal{T}$-modules for all $i, j\in \mathcal{D}_\ell$.

In fact, $P_i^\ell$ is $R$-linearly spanned by $\bigcup _{w\in \{0,\ell\},j\in \mathcal{C}_w}\{b_{ji}^w\}$.
We define a map $$f:P_i^\ell\lra P_{j}^\ell, \ b_{ui}^w\mapsto b_{uj}^w, \forall \, w\in \{0,\ell\},u\in \mathcal{C}_w.$$
Then $f$ is an isomorphism of $\mathcal{T}$-modules, and therefore $P_{i}^\ell\simeq P_{j}^\ell$ as $\mathcal{T}$-modules.

(iii) For $\ell\in [r]$, we fix an element $\ell_0 \in \mathcal{D}_{\ell}$ and define $P_{\ell}$ $:= P^{\ell}_{\ell_0}$. Then $P_{\ell}$ is indecomposable, $P_m\not\simeq P_n$ for $0\le m\neq n\le r$, and $\Hom_{\mathcal{T}}(P_m,P_n)\simeq Rb^0_{m_0,n_0} + \delta_{mn}Rb^m_{m_0, m_0}$, where $\delta_{ij}$ is the Kronecker symbol.

Indeed, for $0\le i,j \le r$, $\Hom_{\mathcal{T}}(P_i,P_j)\simeq b^i_{i_0 i_0}\mathcal{T}b^j_{j_0,j_0}$, that is, $\Hom_{\mathcal{T}}(P_i,P_j)$ is spanned $R$-linearly by $\{b^0_{i_0,j_0},\delta_{ij}b^i_{i_0,i_0}\}$. If $i=j$,  then $\End_{\mathcal{T}}(P_i)$ has a unique nonzero idempotent element $b^i_{i_0,i_0}$.  This implies
that $P_i$ is indecomposable.

Now we show that $P_i\not\simeq P_j$ if $i\neq j$. First, we show that no homomorphism $\varphi\in \Hom_{\mathcal{T}}(P_i,P_j)$ is an isomorphism for $i>j$. Actually, $\Hom_{\mathcal{T}}(P_i,P_j)$ is $1$-dimensional and $\varphi$ is given by the right multiplication by $\lambda b^0_{i_0,j_0}$ with $\lambda\in R$. Since $k_{i_0}=2=p$ and $P_i$ is  spanned $R$-linearly by $\bigcup _{w\in \{0,i\},u\in \mathcal{C}_w}\{b_{u,i_0}^w\}$, we get
$$(b^0_{u i_0})\varphi= b^0_{u i_0}(\lambda b^0_{i_0, j_0})=\lambda k_{i_0} b^0_{u, j_0}= 2 \lambda b^0_{u, j_0}=0$$for $u\in \mathcal{C}_0$. Thus $\varphi$ has a nonzero kernel and is not an isomorphism. This show $P_i\not\simeq P_j$ as $\mathcal{T}$-modules.

If $f: P_i\to P_j$ is an isomorphism of $\mathcal{T}$-modules for $i<j$, then $f^{-1}:P_j\to P_i$ is an isomorphism of $\mathcal{T}$-modules with $j>i$. This contradicts to the foregoing discussion. Thus $P_i\not\simeq P_j$ for all $i\ne j$.

Now, it follows from (i)-(iii) that $\{P_0,P_1,\cdots,P_r\}$ is a complete set of non-isomorphic indecomposable projective $\mathcal{T}$-modules.

(iv) Let $b_\ell:=b^\ell_{\ell_0,\ell_0}$. Then the basic algebra of $\mathcal{T}$ is
\begin{align*}
        \Gamma &:=\End_\mathcal{T}(P_0\oplus P_1\oplus\cdots\oplus P_r)
        =\End_\mathcal{T}(\mathcal{T}b_0\oplus \mathcal{T}b_1\oplus\cdots\oplus \mathcal{T}b_r)
        \simeq
        \begin{pmatrix}
             b_0\mathcal{T}b_0 & \cdots & b_0\mathcal{T}b_r \\
             \vdots & \ddots & \vdots \\
             b_r\mathcal{T}b_0 & \cdots & b_r\mathcal{T}b_r
        \end{pmatrix}\\
        &\simeq
        \begin{pmatrix}
            Rb^0_{0_0,0_0} & Rb_{0_0,1_0}^0 & Rb_{0_0,2_0}^0 & \cdots & Rb_{0_0,r_0}^0 \\
            Rb_{1_0,0_0}^0  & Rb_{1_0,1_0}^0+Rb^1_{1_0,1_0}  & Rb_{1_0,2_0}^0  & \cdots & Rb_{1_0,r_0}^0  \\
            Rb_{2_0,0_0}^0  & Rb_{2_0,1_0}^0  & Rb_{2_0,2_0}^0+Rb^2_{2_0,2_0}  & \cdots & Rb_{2_0,r_0}^0  \\
            \vdots & \vdots & \vdots & \ddots & \vdots \\
            Rb_{R_0,0_0}^0  & Rb_{r_0,1_0}^0  & Rb_{r_0,2_0}^0  & \cdots & Rb_{r_0,r_0}^0+Rb^r_{r_0,r_0}
        \end{pmatrix}_{(r+1)\times (r+1)}
     \end{align*}
For $u,v\in[r]$, let $E^1_{uv}(a)$ denote the $(r+1)\times(r+1)$ matrix with $(u+1,v+1)$-entry $a$ and other entries $0$. Obverse that $\bigcup_{u,v\in[r]}\{E^1_{uv}(b^0_{u_0,v_0})\}\cup \bigcup_{w\in \{1,2,\cdots,r\}} \{E^1_{ww}(b^w_{w_0,w_0})\}$ is an $R$-basis of $\Gamma$ and that $\bigcup_{i\in[r]}\{e_i\}\cup\bigcup_{i\in \{1,2,\cdots,r\}}\{\alpha_i,\beta_i\}\cup\bigcup_{i,j\in \{1,2,\cdots,r\}}\{\beta_i\alpha_i\}$ is an $R$-basis of $\Lambda$.

Let $\psi:\Lambda\to\Gamma$ be the $R$-linear map given by
$$e_i\mapsto E^1_{ii}(b^i_{i_0,i_0}),\ \alpha_{m}\mapsto E^1_{0m}(b^0_{0_0,m_0}),\ \beta_{n}\mapsto E^1_{n0}(b^0_{n_0,0_0}),\ \beta_{u}\alpha_{v}\mapsto E^1_{uv}(b^0_{u_0,v_0}),$$
for $i\in[r]$ and $1\le m,n,u,v\le r$. Since $\psi$ sends the $R$-basis of $\Lambda$ bijectively to the one of $\Gamma$ and preserves multiplication of basis elements, $\psi$ is a homomorphism of algebras. This shows that $\psi$ is an isomorphism of algebras.
$\square$

The basic algebra $\Lambda$ of $\mathcal{T}$ is the dual extension of a star with $r+1$ vertices and $r$ arrows directing to the center of the star. For an algebra $A$ over a field $R$ given by a quiver $Q = (Q_0, Q_1)$ with relations $\{\rho_i\mid i\in I\}$, the \emph{dual extension} of  $A$ (see \cite{x1}) is given be the quiver $(Q_0, Q_1\cup Q'_1)$ with relations $\{\rho_i\mid i\in I\}\cup \{\rho'_i\mid i\in I\} \cup \{\alpha'\beta  \mid \alpha, \beta\in Q_1\}$, where $Q_1'$ is the set of opposite arrows in $Q_1$, that is $Q_1':=\{\alpha':j\to i\mid \alpha: i\to j \mbox{ in } Q_1\}$. It was shown that the global dimension of the dual extension of $A$ is the double of the global dimension of $A$ \cite{x3}.

Theorem \ref{tiqhc}(2) shows that the global dimension of the Terwilliger algebra of a quasi-thin scheme is at most $2$. Thus we re-obtain the quasi-heredity of these algebras by a result of Dlab-Ringel which says that finite-dimensional algebras of global dimension at most $2$ are always quasi-hereditary \cite[Theorem 2]{dr}.

The representation dimension of an Artin algebra $A$ was introduced by Auslander \cite{aus} and defined as follows.
$$\rpd(A):=\inf\{\gd(\End_{A}(A\oplus D(A)\oplus M))\mid M\in A\modcat\}.$$
For the representation dimension of the Terwilliger algebra $\mathcal{T}$, it follows from  \cite[Theorem 3.5]{x2} that $\rpd(\mathcal{T})=\rpd(\Lambda)\le 3$. If $\mathcal{T}$ is semisimple, then $\rpd(\mathcal{T})=0$. Aussme that $\mathcal{T}$ is not semisimple. Then $p=2$ and $r\ge 1$. If $r=1$, then  $\rpd(\mathcal{T})=\rpd(\Lambda)=2$. If $r\ge 2$, then $\Lambda$ is representation-infinite by \cite[Lemma 3.4]{x1}, and therefore $\rpd(\mathcal{T})=\rpd(\Lambda)\ge 3$. In this case, $\rpd(\mathcal{T})=\rpd(\Lambda)=3$.

The dominant dimension of an Artin algebra $A$ is defined to be the minimal nonnegative integer $n$ in a minimal injective resolution
$$ 0\lra {}_AA\lra I_0\lra I_1\lra \cdots \lra I_n\lra \cdots$$of $_AA$, such that $I_n$ is not projective. Related to dominant dimensions, there is a long-standing, famous conjecture, namely the Nakayama conjecture which states that any finite-dimensional algebra over a field with infinite dominant dimension is self-injective (see \cite{N} and \cite[Conjecture (8), p.410]{ars}). This conjecture is still open. However, for the Terwilliger algebras of quasi-thin schemes, the Nakayama conjecture holds true. This can be seen from the following.

\begin{Koro}
    Let $R$ be a field of characteristic $p\ge 0$, and assume that $S$ be a quasi-thin scheme on a finite set $X$, $\mathcal{A}_2$ has $r\ge 0$ equivalence classes. For the Terwilliger $R$-algebra $\mathcal{T}$ of $S$, then the following hold.
$$\dm(\mathcal{T})=\begin{cases} 0, &\mbox { if } p=2 \mbox{ and } r \ge 2, \\ 2, & \mbox { if } p=2 \mbox{ and } r=1,\\ \infty, &\mbox { otherwise}.
\end{cases}$$
\end{Koro}

{\it Proof.} If $\mathcal{T}$ is semisimple, then $\dm(\mathcal{T})=\infty.$ By Proposition \ref{yj1}(3), $\mathcal{T}$ is semisimple if and only if $p\ne 2$ or $p=2$ and $S$ is thin. In the latter, $\mathcal{T}=M_X(R)$. Now let $p=2$ and assume that the quasi-thin scheme $S$ is not thin. In this case, we have $r\ge 1$. By Theorem \ref{tiqhc}(2),  $\dm(\mathcal{T})=
2$ if $r=1$, and $0$ if $r\ge 2$ because in the latter case, no injective $\mathcal{T}$-modules are projective. $\square$

Thus we have the following.

\begin{Koro} The Nakayama conjecture holds true for the class of Terwilliger algebras of quasi-thin schemes over any field.

\end{Koro}

\medskip
As another consequence, we re-obtain the following result in \cite{yj1}.

\begin{Koro} The Jacobson radical of the Terwilliger algebra of a quasi-thin scheme over a field has nilpotent index at most $3$.\end{Koro}

{\it Proof.} Since the statement is true for semisimple Terwilliger algebras of quasi-thin schemes, we have to consider the case of Theorem \ref{tiqhc}(2). In this case, the basic algebras of the Terwilliger algebra of a quasi-thin scheme is radical-cube-zero by Theorem \ref{tiqhc}(2). It is known that Morita equivalent algebras $A$ and $B$ have the isomorphic lattices of left ideals (respectively, ideals). In fact, if a bimodule $_AM_B$ defines a Morita equivalence between $A$ and $B$, then the isomorphism from the lattice of ideals of $A$ to the one of $B$ is given by $I\mapsto I'$ for $I, I'$ ideals in $A$ and $B$, respectively, such that $IM=MI'$ (see \cite[Chapter 2]{bass}). This implies that $\rad(A)^nM=M\rad(B)^n$ for all $n\ge 0$, where $\rad(A)$ stands for the Jacobson radical of $A$. Hence the nilpotent indices of the radicals of $A$ and $B$ are equal. Thus the Jacobson radical of the Terwilliger algebra of a quasi-thin scheme over a field has a nilpotent index at most $3$. $\square$

\medskip
{\bf Acknowledgement.}
The research work was partially supported by the National Natural Science Foundation of China (12031014, 12226314). The authors thank Dr. Yu Jiang for his talk at CNU and discussions on Terwilliger algebras.

 {\footnotesize
 \smallskip
 Zhenxian Chen

School of Mathematical Sciences, Capital Normal University,  100048  Beijing, P. R. China;

 {\tt Email: czx18366459216@163.com}

 \smallskip
 Changchang Xi,

 School of Mathematical Sciences, Capital Normal University, 100048
 Beijing, P. R. China; \&
School of Mathematics and Statistics, Shaanxi Normal University, 710119 Xi'an, P. R. China

 {\tt Email: xicc@cnu.edu.cn}
}
\end{document}